\documentclass[letterpaper, 10 pt, conference]{ieeeconf}  % Comment this line out if you need a4paper
\IEEEoverridecommandlockouts                              % This command is only needed if you want to use the \thanks command
\usepackage{graphicx}
\usepackage{subfigure}
\usepackage{wrapfig} 
\usepackage{url}
\usepackage{hyperref}
\usepackage{indentfirst}
\usepackage{amsmath}
\usepackage{amssymb}
\usepackage{amsfonts}
\usepackage{textcomp}
\usepackage{color}
\usepackage[table]{xcolor}
\usepackage{xcolor}
\usepackage{array}
\usepackage{algorithm} 
\usepackage{algpseudocode}
\usepackage{mathtools}
\DeclarePairedDelimiterX{\inp}[2]{\langle}{\rangle}{#1, #2}
\newcommand{\utwi}[1]{\mbox{\boldmath $#1$}}
\newcommand{\btheta}{{\utwi{\theta}}}
\newcommand{\bvarphi}{{\utwi{\varphi}}}
\newcommand{\norm}[1]{\left\lVert#1\right\rVert}

\newtheorem{theorem}{Theorem}

\title{\LARGE \bf
Data-Driven and Online Estimation of  Linear Sensitivity Distribution Factors: A Low-rank Approach
}

\author{Ana M. Ospina and Emiliano Dall'Anese% <-this % stops a space
\thanks{A. Ospina and E. Dall'Anese are with Department of Electrical, Computer and Energy Engineering, University of Colorado, Boulder, CO, USA, 80309 {\tt\small \{ana.ospina, emiliano.dallanese\}}\tt\small{@colorado.edu}}%
}

\begin{document}

\maketitle
\thispagestyle{empty}
\pagestyle{empty}

%%%%%%%%%%%%%%%%%%%%%%%%%%%%%%%%%%%%%%%%%%%%%%%%%%%%%%%%%
\begin{abstract}
Estimation of sensitivity matrices in electrical transmission systems allows grid operators to evaluate in real-time how changes in power injections reflect into changes in power flows. In this paper, we propose a robust low-rank minimization approach to estimate sensitivity matrices based on measurements of power injections and power flows. An online proximal-gradient method is proposed to estimate sensitivities on-the-fly from real-time measurements. The proposed method obtains meaningful estimates with fewer measurements when the regression model is underdetermined, in contrast with existing methods based on least-squares approaches. In addition, our method can also identify faulty measurements and handle missing data. In this work, convergence results in terms of dynamic regret are presented. Numerical tests corroborate the effectiveness of the novel approach and the robustness of missing measurements and outliers.
\end{abstract}

%%%%%%%%%%%%%%%%%%%%%%%%%%%%%%%%%%%%%%%%%%%%%%%%%
\section{Introduction}

Sensitivity factors play an important role in power systems operations and control. In transmission systems, linear sensitivity distribution factors have traditionally been utilized in power systems analysis -- e.g., contingency analysis, generation re-dispatch, and security assessment \cite{PS_DSE_2019}, just to mention just a few. Injection shift factors (ISFs)~\cite{SensitivityPS_1968,Chen_2014} as well as power transfer distribution factors (PTDFs) allow grid operators to estimate line flows in real-time in response to changes in the (net) power injections. Computation of these sensitivities typically relies on either \emph{model-based} or \emph{measurement-based} approaches. As an example of a model-based method, ISFs and the PTDF matrix for transmission systems are typically computed by leveraging the DC approximation \cite{Sauer_1981}. Model-based approaches require accurate knowledge of the network topology (including line impedances and switchgear states), and are not dependent on specific operating points of the network \cite{Sauer_1981}. Measurement-based methods leverage data obtained from phase measurement units (PMUs) or Supervisory Control and Data Acquisition (SCADA) systems, to obtain estimates of the sensitivity matrix using, e.g., a least-squares approach or alternative estimation criteria \cite{Chen_2014}. Measurement-based methods do not require a knowledge of the topology and impedances, and they do not rely on pseudo-measurements obtained via power flow solutions.

Approaches based on the least-square estimation criterion are effective only if one can collect measurements of the net power injections that are ``sufficiently rich''; that is, measurements that lead to a least-square  that is well conditioned~\cite{Chen_2014, refA3, refA1, refA4}. In principle, the LS may be well conditioned when the perturbations of the net power injections can be properly designed by the grid operator; for example, by adopting the probing techniques of~\cite{bhela2019smart} at some or all of the nodes. However, an underdetermined system may emerge when (i) perturbations may not be performed at a sufficient number of nodes (thus, the power variations are simply due to uncontrollable devices); (ii) changes in the power of uncontrollable loads and generation units located throughout the network may lead to correlated measurements~\cite{abdullah2013probabilistic}; and, (iii) when the power network is operating under dynamic conditions due to fluctuations introduced by intermittent renewable generation~\cite{DhopleNoFuel}, the operator may not have time to collect enough measurements before the operating point of the network changes (and, thus, the sensitivities change). Approaches based on least-squares as \cite{Chen_2013, refA3, refA1, refA5} do not consider the case where the collected measurements have outliers that might lead to unreliable estimates.

To address these challenges, this paper proposes a robust nuclear norm minimization method to estimate sensitivities from measurements. The proposed approach is motivated by our observation that certain classes of sensitivity matrices -- including the PTDFs -- can afford a low-rank approximation. Relative to existing methods based on the least-squares approach, the proposed method: (C1) obtains meaningful estimates of the sensitivity matrices  with a smaller number of measurements and when the regression model is underdetermined (this is particularly important in time-varying conditions and in case of switches in the topology); (C2) leverages sparsity-promoting regularization functions to identify faulty measurements;  and, (C3) utilizes a low-rank approach to handle missing data and measurements collected at different rates. 

To adapt to power networks increasingly operating under dynamic conditions (and, hence, having sensitivity matrices that change rapidly over time), the development of real-time algorithms that can estimate the sensitivity matrix on-the-fly from real-time measurements is presented in this paper. In particular, we propose an online proximal-gradient method~\cite{Online-opti} to solve the nuclear norm minimization problem based on  measurements collected from PMUs and SCADA systems at the second or sub-second level. In line with the broad literature on online optimization, convergence results in terms of dynamic regret~\cite{Hall15,Jadbabaie15} are offered. We point out that the proposed algorithm is markedly different from the competing alternative~\cite{akhriev2020}, and relies on an online proximal-gradient method. Lastly, it is also worth recognizing related works such as \cite{GISF_2017}, where the AC equations are perturbed in order to derive a closed-form expression of so-called ``generalized" ISFs. An approach to estimate dynamic distribution factors is introduced in \cite{DDF_2019}, where reduced-order models are used to derive dynamic ISFs and generator participation factors. An example of online convex optimization in power systems is presented in \cite{dataDrivenLR}, for the specific application of estimating load changes in the grid. 

To the best of our knowledge, we present the first work that explores nuclear norm minimization methods to estimate sensitivity matrices in power systems. Additionally, we develop an online algorithm based on a low-rank model to solve a robust sensitivity estimation problem in the power systems context.

%%%%%%%%%%%%%%%%
\section{Preliminaries}

The proposed approach can be leveraged to estimate various sensitivity coefficients in a power grid. These include, for example,  ISFs~\cite{SensitivityPS_1968,Chen_2014} in transmission systems, and voltage sensitivities (with respect to power injections) in distribution networks~\cite{Rigoni18}. In the following, to clearly and concretely explain the proposed approach, we tailor the exposition to the estimation of the PTDFs matrix in transmission systems\footnote{\textit{Notation:} Upper-case (lower-case) boldface letters will be used for matrices (column vectors), and $(\cdot)^\top$ denotes transposition. For a given column vector $\mathbf{x} \in \mathbb{R}^n$, $\norm{\mathbf{x}} := \sqrt{\mathbf{x}^\top\mathbf{x}}$, and $\norm{\mathbf{x}}_1 := \sum_{i=1}^n | x_i |$. Given a matrix $\mathbf{X} \in \mathbb{R}^{n \times l}$, $\text{vec}(\mathbf{X}) \in \mathbb{R}^p$ denotes the column vectorized $\mathbf{X}$ with its columns stacked in order on top of one other and $p := nl$, $\norm{\mathbf{X}}_* := \sum_{i=1}^r \sigma_i (\mathbf{X})$, where $r$ is the rank of $\mathbf{X}$, and $\sigma_i$ represented the singular values of $\mathbf{X}$. A vector of ones (zeros) is represented by $\mathbf{1}$ ($\mathbf{0}$) with the corresponding dimensions. $\mathcal{O}$ refers to the big O notation; that is, given two positive sequences $\{a_k\}_{k = 0}^\infty$ and $\{b_k\}_{k = 0}^\infty$, we say that $a_k = \mathcal{O}(b_k)$ is $\limsup_{k \rightarrow \infty}  (a_k/b_k) < \infty$.
}.  

\subsection{System Model} \label{System_model}

Let $\mathcal{N}:= \{1, \dots, n\}$ be the set of nodes where generators and/or loads are located, and let $\mathcal{L}:= \{1, \dots,l\}$ be the set of transmission lines. Towards this, let $\Delta p_j \in \mathbb{R}$ represent a change in the net active power injection at node $j \in \mathcal{N}$, around a given point $p_j$; then, the vector capturing the change in the active power flow on the lines in response to the change of power $\Delta p_j$ can be approximated as $\mathbf{h}_j \Delta p_j$,  where $\mathbf{h}_j \in \mathbb{R}^{l}$ represent the sensitivity coefficients~\cite{SensitivityPS_1968,Chen_2014}. Discretize the temporal axis as $\{t_k = k T, k \in \mathbb{N}\}$, with $T$ as a given time interval. Let  $\Delta \mathbf{p}_k := [\Delta p_{1 k}, \Delta p_{2 k}, \dots, \Delta p_{n k}]^\top$ be the vector of net active power changes collected at time instant $t_k$ at the $n$ nodes, and define  the sensitivity matrix as $\mathbf{H}_k := [ \mathbf{h}_{1 k} \; \mathbf{h}_{2 k} \; \dots \; \mathbf{h}_{n k} ] \in \mathbb{R}^{l \times n}$. Then, the vector $\Delta \mathbf{f}_k \in \mathbb{R}^l$ representing the change in the power flow on the lines in the network due to $\Delta \mathbf{p}_k$ can be expressed by~\cite{SensitivityPS_1968}
\begin{equation}
\Delta \mathbf{f}_k = \mathbf{H}_k \Delta \mathbf{p}_k,
\label{Vec-Ma_GSF}
\end{equation}
\noindent
where the entry $i,j$ of $\mathbf{H}_k$ represents the sensitivity injection shift factors~\cite{Chen_2014}. Overall, $\mathbf{H}_k $ can be thought as a proxy for the Jacobian of the map $\mathbf{f} = \mathcal{F}(\mathbf{p})$, which yields flows as a function of power injections, calculated at a given point. By considering $m$ measurements\footnote{Here, we consider measurements taken at times $t_{k-m+1}, \ldots, t_{k}$ for exposition simplicity; however, one may use measurements collected at irregular intervals.} we can define the matrices $ \Delta \mathbf{F}_k = [\Delta \mathbf{f}_{k - m+1} \; \dots \; \Delta \mathbf{f}_k] \in \mathbb{R}^{l \times m}$,  and $ \Delta \mathbf{P}_k = [\Delta \mathbf{p}_{k - m+1} \;  \dots \; \Delta \mathbf{p}_{k}] \in \mathbb{R}^{n \times m}$. Then, the following linear system of equations can be written as
\begin{equation}
\Delta \mathbf{F}_k = \mathbf{H}_k \Delta \mathbf{P}_k .
\label{Matrix_GSF}
\end{equation}
\noindent
Based on~\eqref{Matrix_GSF}, the following subsection will review existing approaches based on the least-squares method as well as model-based approaches.

\subsection{Existing Methods} \label{Existing_models}

\subsubsection{Least-squares estimation}

Assuming that $\Delta \mathbf{P}_k$ is known and measurements (or pseudo-measurements) of $\Delta \mathbf{F}_k$ are available, one possible way to estimate $\mathbf{H}_k$ is via a least-squares criterion. For example, a method similar to  \cite{Chen_2014} can be used, where the injection shift factors for a branch were estimated using  PMU measurements obtained in (near) real-time. In particular, borrowing the approach of~\cite{Chen_2014}, $\mathbf{H}_k$ can be obtained at time $t_k$ by solving: 
\begin{equation}
    \mathbf{H}_{LS,k} \in \arg \min_{\mathbf{H} \in \mathcal{H}} \norm{\Delta \mathbf{F}_k - \mathbf{H} \Delta \mathbf{P}_k}^2_F,
    \label{ISF_LSE}
\end{equation}
where $\norm{\cdot}_F$ denotes the Frobenious norm, and $\mathcal{H}$ is a compact set ensuring that each entry $(i,j)$ of the matrix $\mathbf{H}$ satisfies the constraint $h_{\min} \leq [\mathbf{H}]_{ij} \leq h_{\max}$; that is, $\mathcal{H} = [h_{\min}, h_{\max}]^{l n}$ (in this case, $h_{\min} = -1$ and $h_{\max} = 1$). Alternatively, a weighted least-squares method can be utilized when the noise affecting $\Delta \mathbf{F}_k$ is colored or it is not identically distributed across lines. Notice that in~\eqref{ISF_LSE} there are $lm$ measurements and $ln$ unknowns. With this in mind, existing works such as~\cite{Chen_2014} generally assume that $m \geq n$ and that, the matrix $\Delta \mathbf{P}_k$ has full column rank; with these assumptions, one avoids an underdetermined system and, furthermore,~\eqref{ISF_LSE} has a unique solution. In principle, the matrix $\Delta \mathbf{P}_k$ can have full column rank when the perturbations $\{\Delta p_{j k}\}$ can be properly designed by the grid operator~\cite{bhela2019smart}, or when nodes are perturbed in a round-robin fashion. However, this is impractical in a realistic setting (if not infeasible),  because  the grid operator may not have access to controllable devices at each node of the network; moreover, changes in  the power of uncontrollable loads and generation units located throughout the network contribute to $\{\Delta p_{j k}\}$, and this may lead to correlated measurements 
(therefore, system \eqref{Matrix_GSF} becomes underdetermined). In an underdetermined setting, only a minimum-norm solution would be available using a least-squares criterion, which may provide inaccurate estimates of $\mathbf{H}$ (as corroborated in the numerical results in Section~\ref{Simu}). Before presenting the proposed method, we briefly mention a model-based approach. 

%-------------------------------------------%
%
\subsubsection{Model-based method}

For transmission systems, a widely-used model-based approach to calculate the linear sensitivity distribution factors is based on the DC approximation \cite{Sauer_1981}. In particular, by letting  $\mathbf{B} \in \mathbb{R}^{n \times n}$ represent the matrix of line series susceptances of the transmission system, one can calculate the changes in phase angles $\Delta \btheta$ by using the following relation:
\begin{equation}
    \Delta \mathbf{p}_k = \mathbf{B} \Delta \btheta_k.
    \label{DC_flow_delta}
\end{equation}
Define $\mathbf{X} = \text{diag}(\{ -x_{ab} \}) \in \mathbb{R}^{l \times l}$, where $x_{ab}$ represents the line reactance between node $a \text{ and } b \in \mathbb{N}$, and let $\mathbf{A} \in \mathbb{R}^{l \times n}$ be the branch-bus incidence matrix. Then, using the DC power flow formulation, $\mathbf{f}_k$ can be expressed as the linear relation $\mathbf{f}_k = \mathbf{X}^{-1} \mathbf{A} \btheta_k$~\cite{Wood_Wollwnberg}. If we want to express the active power flow perturbation $\Delta \mathbf{f}_k$ due to a change in the phase angles $\Delta \btheta_k$, we can write $\Delta \mathbf{f}_k$ as $\Delta \mathbf{f}_k =  \mathbf{X}^{-1} \mathbf{A} \Delta \btheta_k$. By replacing this equation and \eqref{DC_flow_delta} in \eqref{Vec-Ma_GSF} we obtain the model-based relation for the sensitivity matrix as: $\mathbf{H} = \mathbf{X}^{-1} \mathbf{A} \mathbf{B}^{-1}$. In order to guarantee that the inverse of $\mathbf{B}$ exists, we require that the DC power flow equations for the nodal power balances are linearly independent. Then, taking the node $1$ as the slack bus, and denoting as $\mathbf{B}_r\in \mathbb{R}^{(n-1) \times (n-1)}$ and $\mathbf{A}_r \in \mathbb{R}^{l \times (n-1)}$ the reduced matrices, the final sensitivity matrix is given by $\mathbf{H} = [ \mathbf{0} \quad \mathbf{X}^{-1} \mathbf{A}_r \mathbf{B}_r^{-1} ]$.

In the DC formulation, the sensitivity matrix factors depends only on the topology of the network, and are invariant to changes in the system operation point, such as load and generation perturbations or failures. In this paper, we target sensitivities based on the AC power flows, which depend on the current operating point, as explained in \cite{power_div}. In the following, a low-rank method will be presented, which does not require knowledge of network topology or reactances.

%------------------------------------%
%  LOW-RANK APPROACH
%-----------------------------------%

\section{Low-Rank Approach}
\label{sec:lowrank}

In this section, we present an approach for the estimation of the  matrix $\mathbf{H}$ with the following features: i) it leverages measurements of $\Delta \mathbf{F}_k$ and $\Delta \mathbf{P}_k$ obtained from PMUs, SCADA, or other similar sources, rather than relying on a network model; ii) it allows for obtaining meaningful estimates of $\mathbf{H}$ even when~\eqref{Matrix_GSF} is underdetermined by leveraging a low-rank approximation of $\mathbf{H}$; iii) when $\Delta \mathbf{P}_k$ is full column rank, it yields an estimation accuracy similar to the least-squares estimator and, iv) can handle missing measurements of flows on some lines (i.e., some entries of $\Delta \mathbf{F}_k$ may be missing). We recall that at each time $k$,  $m$ measurements are utilized; consequently the matrices $\Delta \mathbf{F}_k$ and $\Delta \mathbf{P}_k$ processed at time $k$ are constructed as $ \Delta \mathbf{F}_k = [\Delta \mathbf{f}_{k - m+1} \; \dots \; \Delta \mathbf{f}_k] \in \mathbb{R}^{l \times m}$  and $ \Delta \mathbf{P}_k = [\Delta \mathbf{p}_{k - m+1} \;  \dots \; \Delta \mathbf{p}_{k}] \in \mathbb{R}^{n \times m}$, respectively.

For simplicity of exposition, we first consider the case where measurements are error-free (then, we consider noisy measurements and outliers). Based on the model~\eqref{Matrix_GSF}, the nearly low-rank property of $\mathbf{H}$ motivates us to  consider the following affine rank minimization problem (RMP) \cite{Recht_2010}: 
\begin{align} \label{eq:rank_min}
\underset{\mathbf{H} \in \mathcal{H}}{\min} & \;
\text{rank}(\mathbf{H}) \qquad \text{s.t.} \; \text{vec}(\Delta \mathbf{F}_k) = \mathcal{A}(\mathbf{H}),
\end{align}
\noindent where $\mathcal{H}$ is the convex compact set (in the simplest case, the Cartesian product of box constraints), $\text{vec}(\Delta \mathbf{F}_k) \in \mathbb{R}^p$ where $p:= lm$, denotes the vectorized $\Delta \mathbf{F}_k$, and the linear map $ \mathcal{A} : \mathbb{R}^{l \times n} \rightarrow \mathbb{R}^p $ is defined as: $\mathcal{A}(\mathbf{H}) = \textbf{A}_{\mathbf{P},k} \, \text{vec}(\mathbf{H})$, where $\text{vec}(\mathbf{H}) \in \mathbb{R}^d$, $d := ln$, and $\textbf{A}_{\mathbf{P},k}$ is a matrix of dimensions $p \times d$, appropriately  built using the perturbations $\Delta \mathbf{P}_k$. Specifically, matrix $\mathbf{A}_{\mathbf{P},k}$ is the Kronecker product defined by $\mathbf{A}_{\mathbf{P},k} := \Delta \mathbf{P}_k^\top \otimes \mathbf{I}$, where $\mathbf{I}$ is the identity matrix of dimensions $l \times l$. Unfortunately, the rank criterion in~\eqref{eq:rank_min} is in general NP-hard to optimize; nevertheless, drawing an analogy from compressed sensing to rank minimization, the following convex relaxation of RMP~\eqref{eq:rank_min} can be utilized~\cite{Recht_2010}:
\begin{align} \label{eq:nuclear_norm_min}
\underset{\mathbf{H}\in \mathcal{H}}{\min} & \; \norm{\mathbf{H}}_* \qquad  \text{s.t.} \; \text{vec}(\Delta \mathbf{F}_k) = \textbf{A}_{\mathbf{P},k} \, \text{vec}(\mathbf{H}),
\end{align}
where $\norm{\mathbf{H}}_* := \sum_i \sigma_i(\mathbf{H})$ is the nuclear norm of $\mathbf{H}$, with $\sigma_i(\mathbf{H})$ denoting the $i$th singular value of $\mathbf{H}$.  Interestingly, it was shown in \cite{Recht_2010} that, if the constraints of \eqref{eq:nuclear_norm_min} are defined by a linear transformation that satisfies a restricted isometry property condition, the minimum rank solution can be recovered by the minimization of the nuclear norm over the linear space; see the necessary and sufficient condition in \cite{Recht_2010}. Unfortunately, designing the perturbation matrix $\Delta \mathbf{P}_k$ to satisfy the restricted isometry property may not be possible due to uncontrollable loads.

The proposed methodology leverages the relaxation~\eqref{eq:nuclear_norm_min} to estimate the matrix $\mathbf{H}$ from measurements of $\Delta \mathbf{F}_k$ induced by the perturbations in the net power injections $\Delta \mathbf{P}_k$.  Assuming that the measurements $\Delta \mathbf{F}_k$ are affected by a zero-mean Gaussian noise, a pertinent relaxation of~\eqref{eq:nuclear_norm_min} amounts to the following convex program~\cite{candes2010matrix}:  
\begin{equation}
\underset{\mathbf{H}\in \mathcal{H}}{\min}
 \norm{\text{vec}(\Delta \mathbf{F}_k) - \textbf{A}_{\mathbf{P},k} \, \text{vec}(\mathbf{H})}^2_2 + \lambda \norm{\mathbf{H}}_*,
\label{noise_nuclear}
\end{equation}
\noindent
where $\lambda > 0$ is a given regularization parameter that is used to promote sparsity in the singular values of $\mathbf{H}$ (and, hence, to obtain a low-rank matrix $\mathbf{H}$).

% -------------------------------------------
\subsection{Robustness to outliers}

We further consider the case where some measurements of $\Delta \mathbf{F}_k$ may be corrupted by outliers. This can be due to, for example, faulty readings of PMUs, communication errors, or malicious attacks. To this end, we augment the model~\eqref{Matrix_GSF} as $\Delta \mathbf{F}_k = \mathbf{H}_k \Delta {\mathbf{P}_k} + \mathbf{O}_k + \mathbf{E}_k$, where $\mathbf{E}_k$ is a matrix containing (small) measurement errors and $\mathbf{O}_k$ is a matrix containing measurement outliers~\cite{zhou2010stable,Mateos_sparse_2013}. When no outliers are present, $\mathbf{O}_k$ is a matrix with all zeros. Based on this augmented model, estimates of $\mathbf{H}_k$ and $\mathbf{O}_k$ can be sought by solving the following convex problem~\cite{zhou2010stable,Mateos_sparse_2013}:    
\begin{multline}
\underset{\mathbf{H} \in \mathcal{H}, \mathbf{O} \in \mathcal{M}}{\min} \,
\norm{\text{vec}(\Delta \mathbf{F}_k) - \textbf{A}_{\mathbf{P},k} \, \text{vec}(\mathbf{H}) - \text{vec}(\mathbf{O})}^2_2 \\ + \lambda \norm{ \mathbf{H}}_* + \gamma \norm{\text{vec}(\mathbf{O})}_1,
\label{eq:noise_nuclear_outliers}
\end{multline}
where $\norm{\text{vec}(\mathbf{O})}_1 = \sum_{i} |[\text{vec}(\mathbf{O})]_i|$ is the  $\ell_1$-norm of the vector $\text{vec}(\mathbf{O})$,  $\gamma > 0$ is a sparsity-promoting coefficient, and $\mathcal{M}$ are box constraints of the form $o_{\min} \leq [\mathbf{O}_{k}]_{ij} \leq o_{\max}$. Notice that the $\ell_1$-norm is the closest convex surrogate to the cardinality function. Once \eqref{eq:noise_nuclear_outliers} is solved,  the locations of nonzero entries in $\mathbf{O}$ reveal outliers across both lines and time; on the other hand,
the amplitudes quantify the magnitude of the anomalous measurement. It is important to notice that the  parameters $\lambda$ and $\gamma$ control the trade-off between fitting error, rank of $\mathbf{H}$,
and sparsity level of $\mathbf{O}$; in particular, when an estimate of the variance of the measurement noise is available, one can follow guidelines for selecting $\lambda$ and $\gamma$ similar to the ones proposed in~\cite{zhou2010stable}.

\subsection{Missing measurements}

It is worth pointing out that the proposed methodology applies to the case where some measurements in $\Delta \mathbf{F}_k$ in~\eqref{noise_nuclear} are \emph{missing}. This may be due to communication failures or because measurements are collected at different rates. For the latter, one can take the highest measurement frequency (i.e., the $T$ shortest inter-arrival time) as a reference frame, and treat measurements that are received less frequently (i.e., with a larger inter-arrival time) as missing entries. 

In this case, missing measurements are discarded from the least-squares term in~\eqref{eq:noise_nuclear_outliers}~\cite{Mateos_sparse_2013}. In particular, let $\Omega_k \subseteq \{1, 2, \ldots, p\}$ be a set indicating which measurements in the vector $\text{vec}(\Delta \mathbf{F}_k)$ are available at time $t_k$; for example, if the measurement for the line $1$ is missing at time $t_{k-m+1}$ and $t_{k-m+2}$, then, $\Omega_k \subseteq \{2, 3, \ldots, l, l+2, \ldots, p\}$. Let $\mathcal{P}_{\Omega_k}$ be a time-varying vector sampling operator, which sets the entries of its vector argument not indexed by $\Omega_k$ to zero and leaves the other entries unchanged. Then,~\eqref{eq:noise_nuclear_outliers} can be written as:  
\begin{multline}
\underset{\mathbf{H} \in \mathcal{H}, \mathbf{O} \in \mathcal{M}}{\min} \,
\norm{\mathcal{P}_{\Omega_k}\left\{\text{vec}(\Delta \mathbf{F}_k) - \textbf{A}_{\mathbf{P},k} \, \text{vec}(\mathbf{H}) - \text{vec}(\mathbf{O}) \right\}}^2_2 \\ + \lambda \norm{ \mathbf{H}}_* + \gamma \norm{\text{vec}(\mathbf{O})}_1,
\label{eq:noise_nuclear_outliers_2}
\end{multline}
where, of course, missing measurements are not accounted for in the least-squares term.

%--------------------------------------------------

\section{Data-Driven Online Estimation} \label{Online}

Based on  $\Delta \mathbf{P}_k$ and $\Delta \mathbf{F}_k$, which collect measurements of new power injections and power flows acquired at time steps $k-m+1, \ldots, k$, an estimate of $\mathbf{H}_k$ can be obtained by solving the convex problem~\eqref{eq:noise_nuclear_outliers} using existing batch solvers for non-smooth convex optimization problems. When the power network is operating under dynamic conditions, for example, due to swings in the net power due to intermittent renewable generation and uncontrollable loads~\cite{DhopleNoFuel}, the sensitivity matrix $\mathbf{H}_k$ may rapidly change over time (since, in general, it depends on the current operating points~\cite{Chen_2014}, \cite{ power_div}); in these dynamic conditions, it may not be possible to solve~\eqref{eq:noise_nuclear_outliers} sufficiently fast due to underlying computational complexity considerations, and a solution of~\eqref{eq:noise_nuclear_outliers} generated by batch solvers can be outdated. That is, by the time the solution is produced, the operating conditions of the network (and, hence, $\mathbf{H}_k$) have changed. This  aspect  motivates the development of an online algorithm that estimates $\mathbf{H}_k$ based on streams of measurements and identifies outliers ``on the fly,'' as explained in this section. 

Measurements are assumed to arrive at times $\{t_k = k T, k \in \mathbb{N}\}$, with $T$ the inter-arrival time (e.g, $T$ could be one second or a few seconds~\cite{DhopleNoFuel}); suppose further that measurements are processed over a sliding window $\mathcal{T}_k = \{t_{k-m+1}, \ldots, t_{k}\}$. Then, at each instant $t_k$, the matrix $\mathbf{H}_k$ can be estimated via~\eqref{eq:noise_nuclear_outliers}, which is re-written here as the following time-varying problem~\cite{Online-opti}:  
\begin{subequations}
\label{eq:timevarying_opt}
\begin{equation}
    (\mathbf{H}_k^*, \mathbf{O}_k^*) \in \arg \underset{\mathbf{H}_k \in \mathbb{R}^{l \times n}, \mathbf{O}_k \in \mathbb{R}^{l \times m}} {\text{min}} \; f_k(\mathbf{H}_k, \mathbf{O}_k) , \hspace{.4cm} \forall \, k T
\end{equation}
where $f_k(\mathbf{H}_k, \mathbf{O}_k) := s_k(\mathbf{H}_k, \mathbf{O}_k)
  + g_k(\mathbf{H}_k, \mathbf{O}_k)$, 
\begin{equation}
  s_k(\mathbf{H}_k, \mathbf{O}_k) := \norm{{\Delta\mathbf{F}}_k - \mathbf{H}_k {\Delta \mathbf{P}}_k - \mathbf{O}_k}_F^2,
\end{equation}
\begin{align}
g_k(\mathbf{H}_k, \mathbf{O}_k) : = & \lambda_k  \norm{\mathbf{H}_k}_* + \gamma_k \norm{\text{vec}(\mathbf{O}_k)}_1 \nonumber \\
& + \iota_{\mathcal{H}}(\mathbf{H}_k) + \iota_{\mathcal{M}}(\mathbf{O}_k), \label{eq:g_k}
\end{align} 
\end{subequations}
with $\iota_{\mathcal{H}}(\mathbf{H})$ the set indicator function for the compact set $\mathcal{H}$ and
 $\iota_{\mathcal{M}}(\mathbf{O})$ the set indicator function for the compact set $\mathcal{M}$.  The goal posed here is to develop an online algorithm that can track a solution $\{\mathbf{H}_k^*, \mathbf{O}_k^*\}_{k \in \mathbb{N}}$ and the trajectory of optimal value functions $\{f_k^* := f_k(\mathbf{H}_k^*, \mathbf{O}_k^*)\}_{k \in \mathbb{N}}$ by processing measurements in a  sliding window fashion. In the following, let  $\mathbf{o}_k = \text{vec}(\mathbf{O}_k)$, and $\mathbf{x}_k = [\text{vec}(\mathbf{H}_k)^\top, \mathbf{o}_k^\top]^\top  \in \mathcal{X}_k : = \mathcal{H} \times \mathcal{M}$ for brevity.  Notice that $s_k(\mathbf{x}_k)$ is closed, convex and proper, with a $L_k$-Lipschitz continuous gradient at each time $t_k$; on the other hand, $g_k(\mathbf{x}_k)$ is a lower semi-continuous proper convex function. Lastly, the function attains a finite minimum at a certain $\mathbf{x}_k^*$. Given this particular structure of~\eqref{eq:timevarying_opt}, we propose to use an online  proximal-gradient  algorithm~\cite{Online-opti} to solve~\eqref{eq:timevarying_opt} under streams of measurements. Assuming that, because of communication delays and computational considerations, one step of the  algorithm can be performed within an interval $T$ (which coincides with the inter-arrival rate of the measurements), the online  proximal-gradient  algorithm amounts to the sequential execution of the following step:
\begin{align} \label{prox_update}
     \mathbf{y}_k &=  \mathbf{x}_{k-1} - \alpha \nabla_\mathbf{x}  s_k(\mathbf{x}_{k-1}), \quad
     \mathbf{x}_{k} = \text{prox}_{g_k,\mathcal{X}}^{\alpha} \{ \mathbf{y}_k \},
\end{align}
 where $\alpha > 0 $ is the step size, and 
the proximal operator is defined over the non-differentiable function $g_k$ as \cite{Beck_prox}
\begin{equation}
    \text{prox}_{g}^{\alpha} \{ \mathbf{y} \} := \underset{\mathbf{x}} {\text{arg min}} \Big \{g(\mathbf{x}) + \frac{1}{2 \alpha} \norm{\mathbf{x} - \mathbf{y}}^2  \Big \} . 
\end{equation}
Notice that, if we re-write the function $s_k$ as
\begin{equation}
    s_k(\mathbf{x}_k) = \norm{\Delta \mathbf{f}_k - \mathbf{A_{Ps}}_{,k} \mathbf{x}_k}^2,
\end{equation}
where $\Delta \mathbf{f}_k = \text{vec}(\Delta \mathbf{F}_k)$, $\mathbf{A_{Ps}}_{,k} = [\mathbf{A_{P}}_{,k}, \mathbf{I}]$, and $\mathbf{x}_k$ defined as before, then,  $\nabla_\mathbf{x} s_k$ is given by
\begin{equation}
    \nabla_\mathbf{x} s_k(\mathbf{x}_k) = 2 \mathbf{A_{Ps}}_{,k}^\top \left( \mathbf{A_{Ps}}_{,k
    } \mathbf{x}_k - \Delta \mathbf{f}_k  \right).
    \label{closed_nabla_s}
\end{equation}
Notice that the proximal operator in~\eqref{prox_update} is separable across the two variables of interest $\mathbf{H}_k$ and $\mathbf{o}_k$, since the costs associated with the $\mathbf{H}_k$ and $\mathbf{o}_k$ are decoupled; therefore, the proximal mappings for $\mathbf{H}_k$ and $\mathbf{o}_k$ can be computed separately. In particular, one has that, 
\begin{align}
     \mathbf{H}_{k} & = \text{prox}_{\lambda_k   \norm{\cdot}_* + \iota_{\mathcal{H}} }\{\mathbf{Y}_{H,k}\},  \label{prox_update_H}\\
     \mathbf{o}_{k} & = \text{prox}_{\gamma_k   \norm{\cdot}_1 + \iota_{\mathcal{M}} }\{\mathbf{y}_{o,k}\},
\end{align}
with $\mathbf{Y}_{H,k}$ and $\mathbf{y}_{o,k}$ extracted from the stacked vector $\mathbf{y}_k$ in~\eqref{prox_update}, then, $\mathbf{o}_{k}$ admits a closed-form solution, given by:
\begin{equation}
    \mathbf{o}_{k} = [\mathcal{S}_\gamma (\mathbf{y}_{o,k})]_{o_{\min}}^{o_{\max}},
    \label{prox_norm1}
\end{equation}
where $[x]_a^b = \max\{\min\{x,b\},a\}$, and the thresholding operator $\mathcal{S}_\gamma$ is defined as:
\begin{align}
    \mathcal{S}_\gamma (\mathbf{y}) & = \max\{|\mathbf{y}| - \gamma \mathbf{1}, \mathbf{0}\} \odot \text{sgn}(\mathbf{y}) \nonumber \\
    & = \Bigg\{ \begin{matrix} \mathbf{y} - \gamma \mathbf{1}, 
    & \textrm{~if~} \mathbf{y} \geq \gamma \mathbf{1}, \\ 
    0, & \textrm{\,~~if~} |\mathbf{y}|< \gamma \mathbf{1}, \\ \mathbf{y} + \gamma \mathbf{1}, & \textrm{~~~if~} \mathbf{y} \leq -\gamma \mathbf{1}. \end{matrix}
\end{align}

With the previous definitions in place, the online proximal-gradient algorithm for the robust estimation of the sensitivity matrix is tabulated as Algorithm \ref{algo}. %We stress that, by using Algorithm \ref{algo}, we can estimate the sensitivity matrix robustly over a sliding window $\mathcal{T}_k$, this adapting to changing operational points of the power system.  
\begin{algorithm}
    \caption{Online robust estimation of sensitivity matrices}\label{algo}
    \begin{algorithmic}
        \For {$k = m, m+1,\dots, $}
        \State \textbf{[S1]} Collect $\Delta \mathbf{f}_k$ and  $\Delta \mathbf{p}_k$
        \State \textbf{[S2]} Build $\Delta \mathbf{F}_k$ and  $\Delta \mathbf{P}_k$ based on $\{\Delta \mathbf{f}_k, \Delta \mathbf{p}_k\}_{k \in \mathcal{T}_k}$ 
        \State \textbf{[S3]} Compute $\mathbf{y}_k$ via~\eqref{prox_update}
        \State \textbf{[S4]} Update $\mathbf{H}_{k}$ via~\eqref{prox_update_H} 
        \State \textbf{[S5]} Update $\mathbf{o}_{k}$ via~\eqref{prox_norm1} 
		\State Go to \textbf{[S1]}
        \EndFor
    \end{algorithmic} 
\end{algorithm}

In order to analyze the estimation accuracy of Algorithm~\ref{algo} performance, the dynamic regret metric is considered here; see, e.g. \cite{Hall15,Jadbabaie15,Online-opti}. In particular, it is defined as: 
\begin{equation*}
    \text{Reg}_k := \sum_{i=1}^k \left[ f_i(\mathbf{x}_i) - f_i(\mathbf{x}_i^*) \right],
\end{equation*}
where we recall $f_i$ is the cost function in~\eqref{eq:noise_nuclear_outliers} (see also~\eqref{eq:timevarying_opt}).
The dynamic regret is an appropriate performance metric for time-varying problems with a cost that is convex, but not necessarily strongly convex~\cite{Online-opti}.  To derive bounds on the dynamic regret, it is first necessary to introduce a ``measure'' of the temporal variability of \eqref{eq:timevarying_opt}. One possible measure is:
\begin{equation}
    \omega_k := \norm{\mathbf{x}_k^*-\mathbf{x}_{k-1}^*},
    \label{eq:omega}
\end{equation}
along with the so-called ``path length'': 

\begin{equation}
    \Omega_k := \sum_{i=1}^k \omega_i, \qquad  \Bar{\Omega}_k := \sum_{i=1}^k \omega_i^2. 
    \label{eq:path_length}
\end{equation}
Recall that the least-squares term $s_k(\mathbf{x}_k)$ is closed, convex and proper, with a $L_k$-Lipschitz continuous gradient at each time $t_k$, and that $g_k(\mathbf{x}_k)$ is a lower semi-continuous proper convex function. Then, by using the definitions~\eqref{eq:omega}--\eqref{eq:path_length} and leveraging bounding techniques  similar to~\cite{Inexact_prox_online}, the following result can be obtained.  

\begin{theorem} 
\label{thm:convergence}
Suppose that the step size $\alpha$ is chosen such that $\alpha \leq 1/L$, with $L := \max\{L_k\}$. Then, the dynamic regret of Algorithm~\ref{algo} has the following limiting behavior:
\begin{equation}
   \frac{1}{k} \text{Reg}_k = \mathcal{O}(1  +  k^{-1} \Omega_k + k^{-1} \bar{\Omega}_k).
\end{equation}
\hfill $\Box$
\end{theorem} 
\emph{Proof.} See the Appendix.

\vspace{.1cm} 

\noindent Note that: 

\noindent $\bullet$ When the sensitivity matrix changes over time, $\Omega_k$ and $\bar{\Omega}_k$ grow as $\mathcal{O}(k)$. Therefore,  $(1/k)\text{Reg}_k = \mathcal{O}(1)$; that is, the sensitivity matrix can be estimated within a bounded error even in the considered online setting~\cite{Online-opti}. 

\noindent $\bullet$ A no-regret result (i.e., $(1/k) \text{Reg}_k$ asymptotically goes to $0$) can not be obtained in general.  

\noindent $\bullet$ If the sensitivity matrix is constant, then one trivially has that  $\text{Reg}_k$ approaches $0$ asymptotically, thus recovering convergence results for the batch proximal-gradient method.

%\begin{remark}
%In the paper, we assume that that $\Delta \mathbf{P}_k$ is either known or can be measured with negligible noise; on the other hand, $\Delta \mathbf{F}_k$ is noisy and may contain outliers. In principle, the proposed approach could be extended to handle noise and outliers in $\Delta \mathbf{P}_k$ by replacing the least-squares term with a total least-squares (TLS) criterion; see, for example,~\cite{markovsky2007overview}. However, the resultant cost function is in this case nonconvex (because of bilinear terms); the challenges rely on the model for the  trajectories of the critical points of the cost function.
%\end{remark}

%-------------------------------------
% SIMULATION
%------------------------------------
%
\section{Simulation Results} \label{Simu}
In this section, the estimation accuracy of the proposed methodology is assessed, for both batch and online implementation.  The following transmission networks are considered: i) Western Electricity Coordinating Council (WECC) 3-machine 9-bus transmission system~\cite{matpower_2019}, and ii) the synthetic South Carolina 500-bus transmission power system model~\cite{Birchfield}. In both cases, MATPOWER is used to compute the power-flow solutions~\cite{ matpower_2019}.

\subsection{Batch estimation} \label{transmission_ex}

Fluctuations in active power injections around a given operating point are simulated as in  \cite{Chen_2014}. In particular, the injection at node $j$, denoted by $p_j$, is given by $p_j[k] = p_j^0[k] + \sigma_{N1} p_j^0[k] \eta_1 + \sigma_{N2} \eta_2$, where $p_j^0[k]$ is the nominal power injection at node $j$ at instant $k$, and ($\eta_1$, $\eta_2$) are random values, where $\eta_1 \sim \mathcal{N}(0, \sigma_{N1})$ and $\eta_2 \sim \mathcal{N}(0, \sigma_{N2})$ for standard deviations $\sigma_{N1} = \sigma_{N2} = 0.1$; see  \cite{Chen_2014} for details. For each time $k$, we take the difference between consecutive line flow measurements to obtain $\Delta \mathbf{F}_k$. We obtained $\Delta \mathbf{P}_k$ by taking the differences between consecutive values of active power injections in each node. The batch optimization problems~\eqref{noise_nuclear} and~\eqref{eq:noise_nuclear_outliers} can also be solved efficiently using the proximal-gradient method (i.e., a batch version of Algorithm~\ref{algo}); see, for example,~\cite{cevher2014convex} (and references therein) for standard computational times of proximal-gradient methods.

\begin{figure}[!t]
\centering
  \subfigure[]{\includegraphics[width=0.41\textwidth]{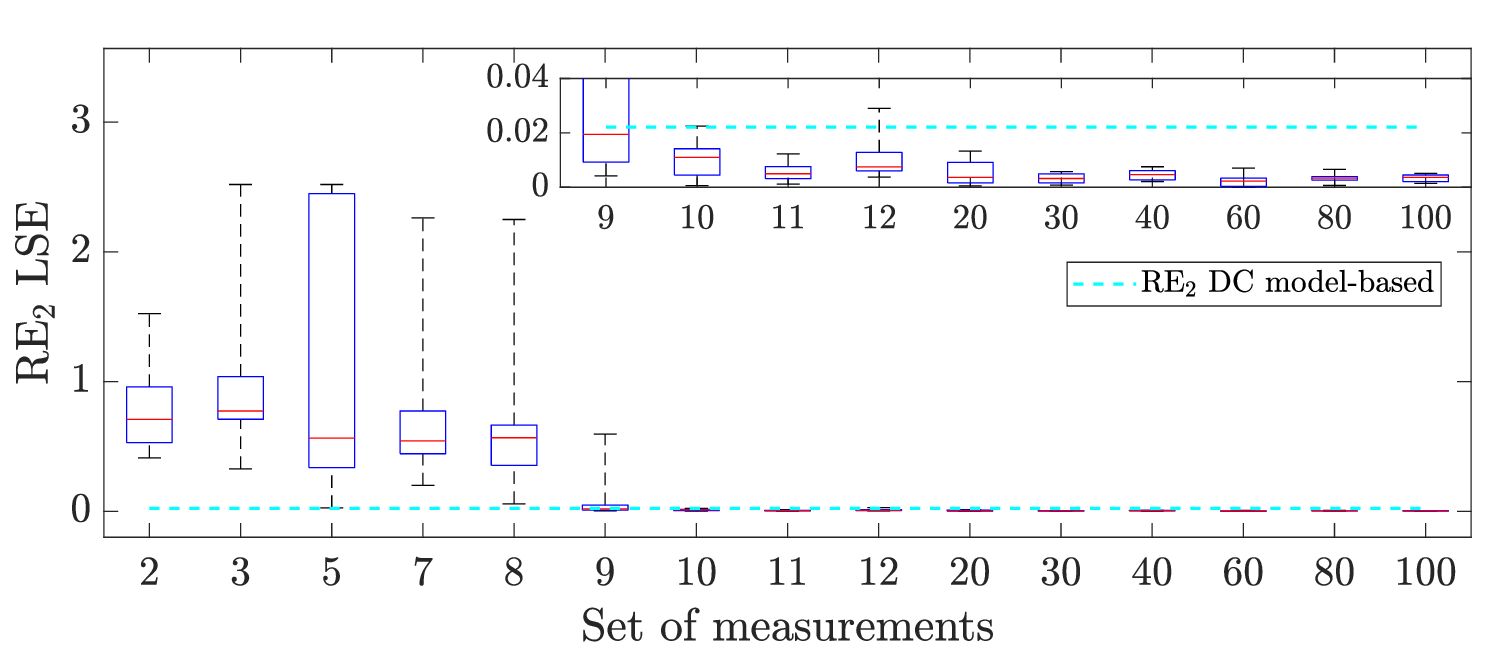}}
  \subfigure[]{\includegraphics[width=0.41\textwidth]{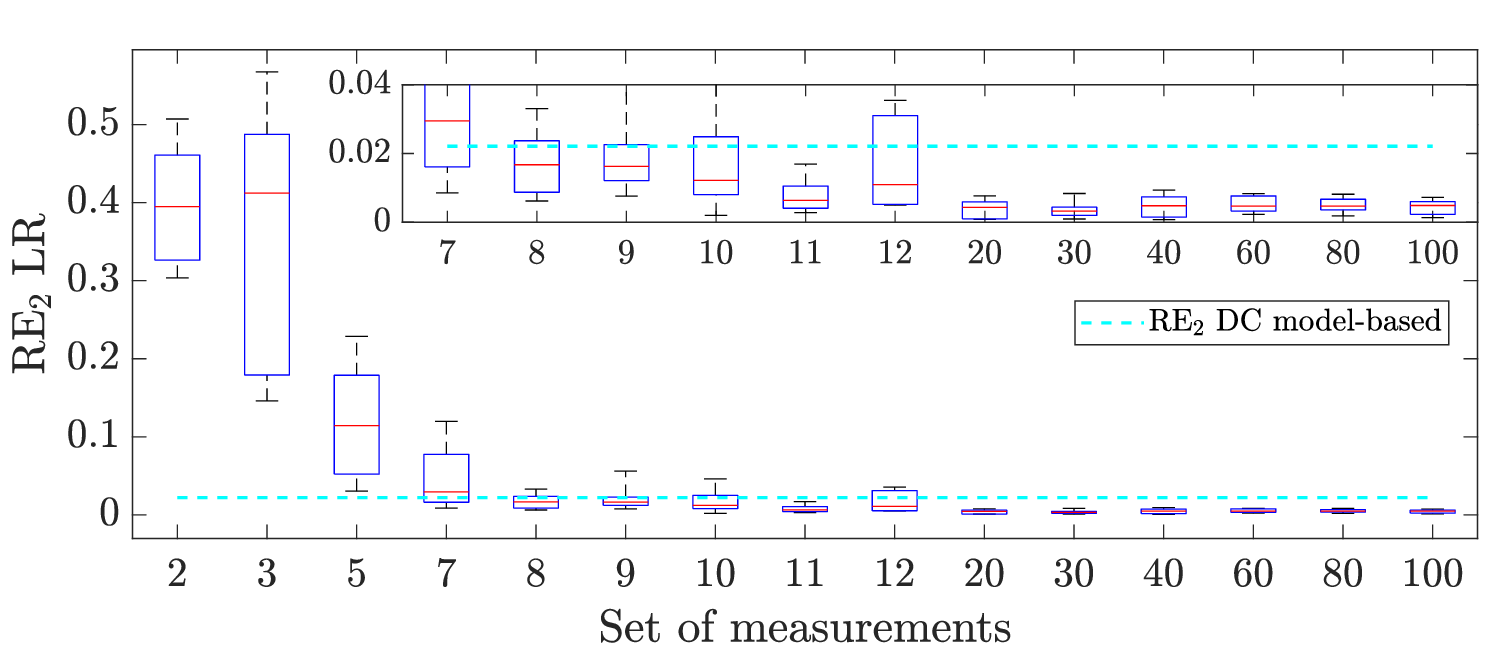}}
  \vspace{-.3cm}
  \caption{Case - 9-bus: Box plot of the relative error (RE) over 10 trials for the estimation of the sensitivity matrix under a different number of measurements: (a) RE for the least square estimator and (b) RE for the low-rank approach.}
  \vspace{-.25cm}
  \label{fig:Transmision_9bus}
\end{figure}

\begin{figure}[!t]
\centering
  \subfigure[]{\includegraphics[width=0.41\textwidth]{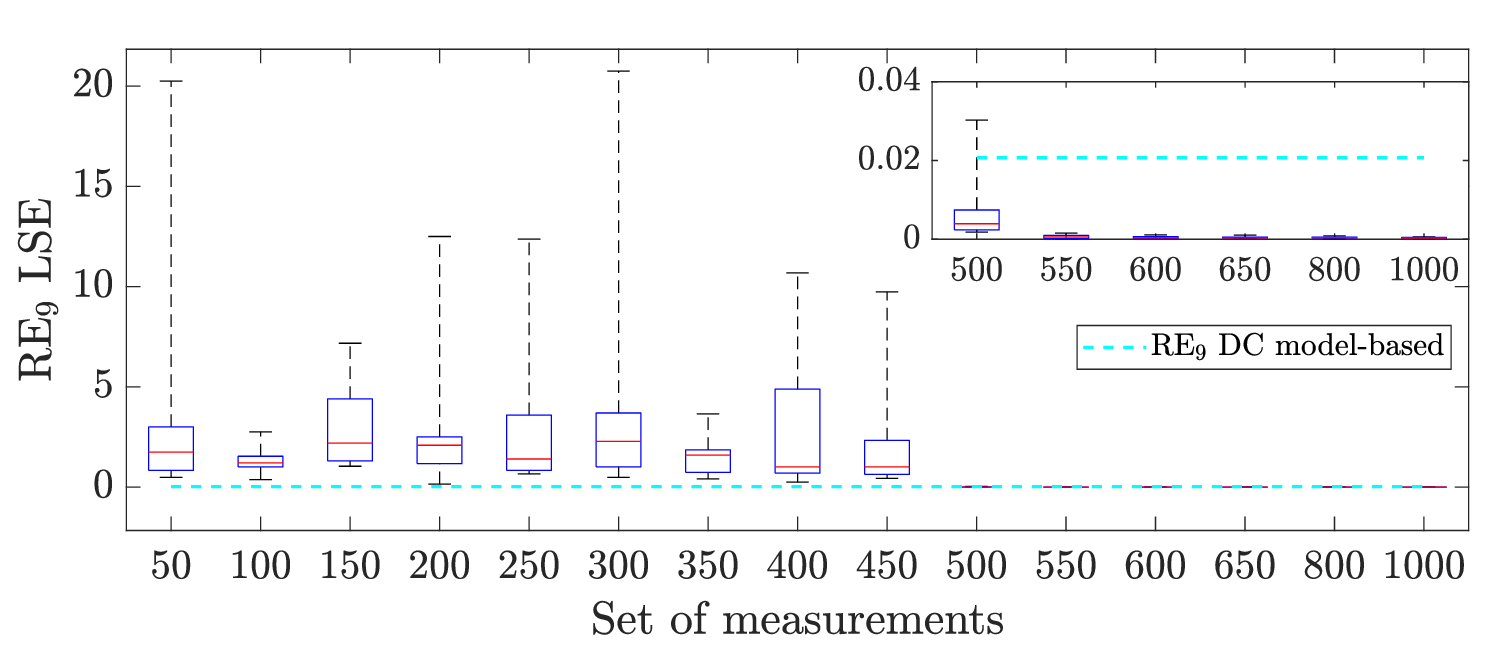}}
  \subfigure[]{\includegraphics[width=0.41\textwidth]{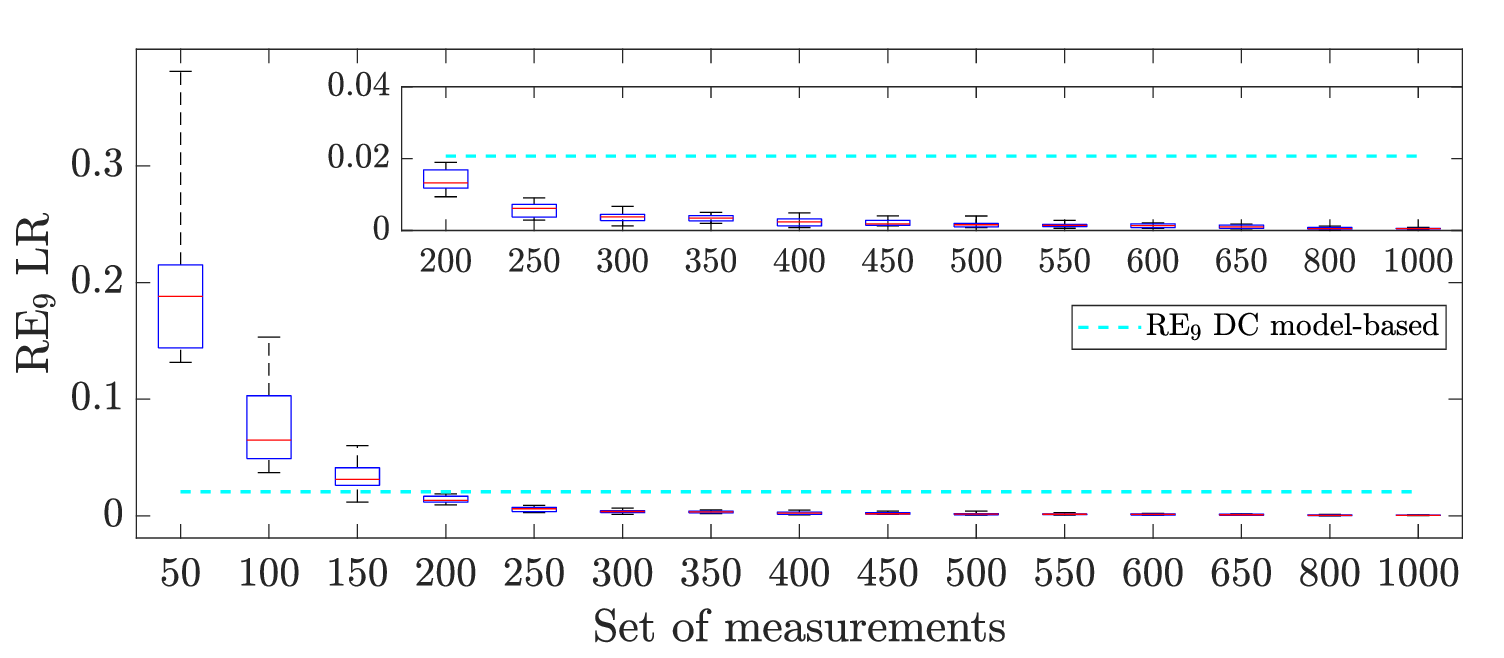}}
  \vspace{-.3cm}
  \caption{Case - 500-bus: Box plot of the relative error (RE) over 10 trials for the estimation of the sensitivity matrix under a different number of measurements: (a) RE for the least-square estimator and (b) RE for the low-rank approach.}
  \vspace{-.25cm}
  \label{fig:Transmision_500bus}
\end{figure}

The performance of the proposed low-rank based approach is considered for both transmission networks. First, the batch method~\eqref{noise_nuclear} is evaluated, when a decrease in generation occurs at generator 2 for the 9-bus, and in generator 9 for the 500-bus transmission system.  Figure \ref{fig:Transmision_9bus} and Figure \ref{fig:Transmision_500bus} compare the performance of the proposed method with the least-squares approach \cite{Chen_2014}. In this case, 10 trials were used, and the relative error (RE) with respect to bus $i$ is defined by $\text{RE}_i = \frac{\norm{\mathbf{h}_{ik} -\mathbf{h}_{ik}^*}}{\norm{\mathbf{h}_{ik}^*}}$ (we use this definition to be consistent with~\cite{Chen_2014}),
where the actual sensitivity of the lines due to the change of generation in bus $i$ is denoted $\mathbf{h}_{ik}^*$, and $\mathbf{h}_{ik}$ specifies the column of the estimated sensitivity matrix $\mathbf{H}_k$ for the bus $i$ obtained form the DC model-based, least-squares, or low-rank approaches. Figure \ref{fig:Transmision_9bus} shows that when the set of measurements is less than 9 (the total number of nodes in this case) the least-squares approach does not give an accurate estimation, because it is underdetermined. In the case of the proposed low-rank method, the median of the relative errors can be just 3\%  even when we have only 6-7 sets of measurements; the proposed method performs better than the model-based approach via DC approximation once we collect 8 measurements. When more than 9 measurements are collected, the proposed method and the least-squares approach have similar performance as expected. Figure~\ref{fig:Transmision_500bus} shows similar behavior for the synthetic South Carolina 500-bus transmission systems, where our method is able to estimate the sensitivity matrix from 200 sets of measurements. In both cases, it is evident that the proposed approach provides  accurate results with fewer measurements than the least-squares approach.

\begin{figure}[!ht]
\centering
  \subfigure[]{\includegraphics[width=0.41\textwidth]{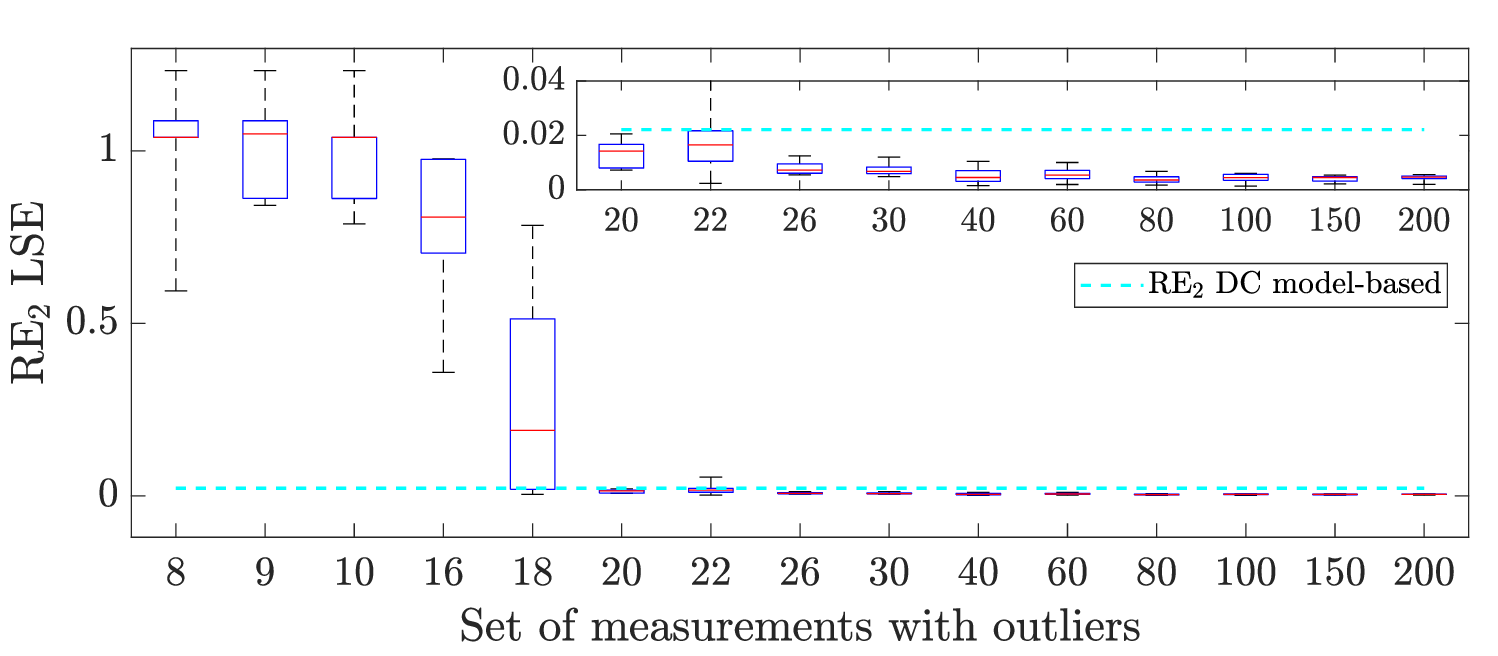}}
  \subfigure[]{\includegraphics[width=0.41\textwidth]{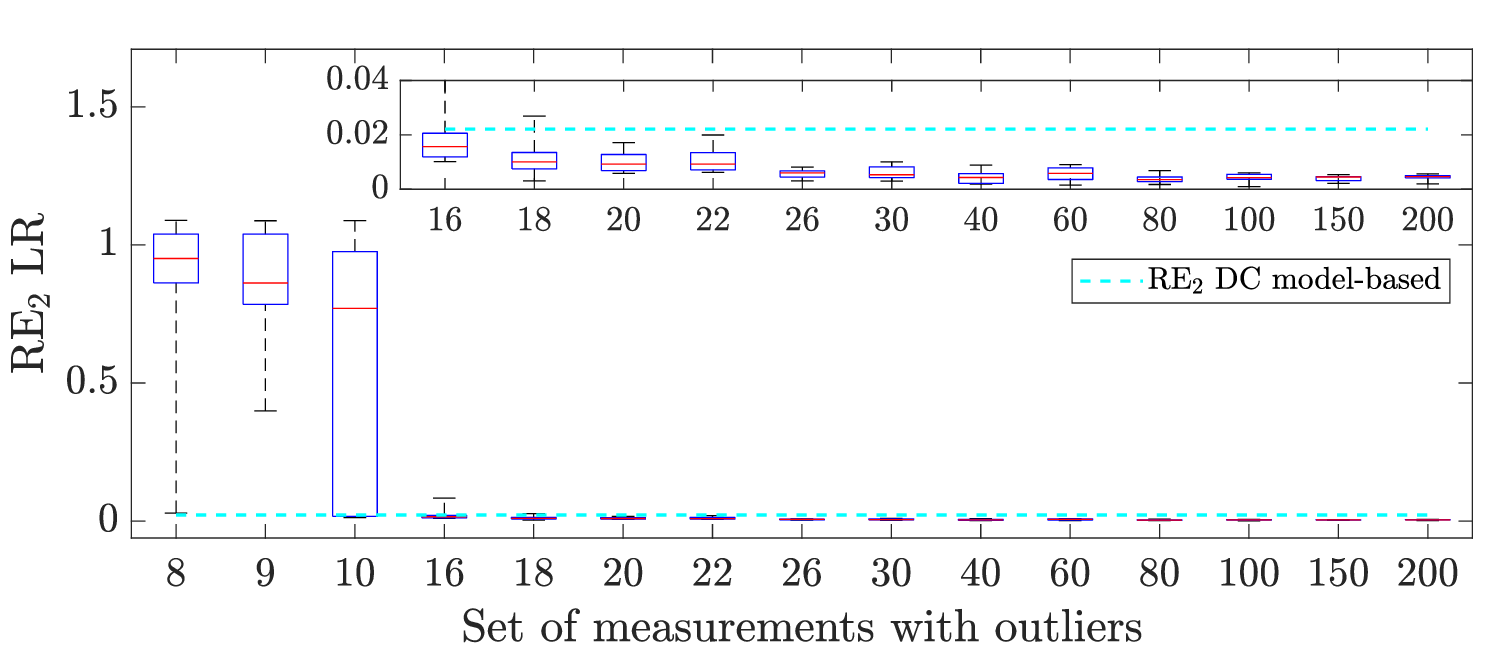}}
  \vspace{-.3cm}
  \caption{Case - 9-bus: Box plot of the relative error (RE)  over 10 trials for the estimation of the sensitivity matrix with outliers under a different number of measurements: (a) RE for the least-square estimator and (b) RE for the low-rank approach.}
  \vspace{-.25cm}
  \label{fig:Figure_Tra_outliers}
\end{figure}

\begin{figure}[!ht]
\centering
  \subfigure[]{\includegraphics[width=0.4\textwidth]{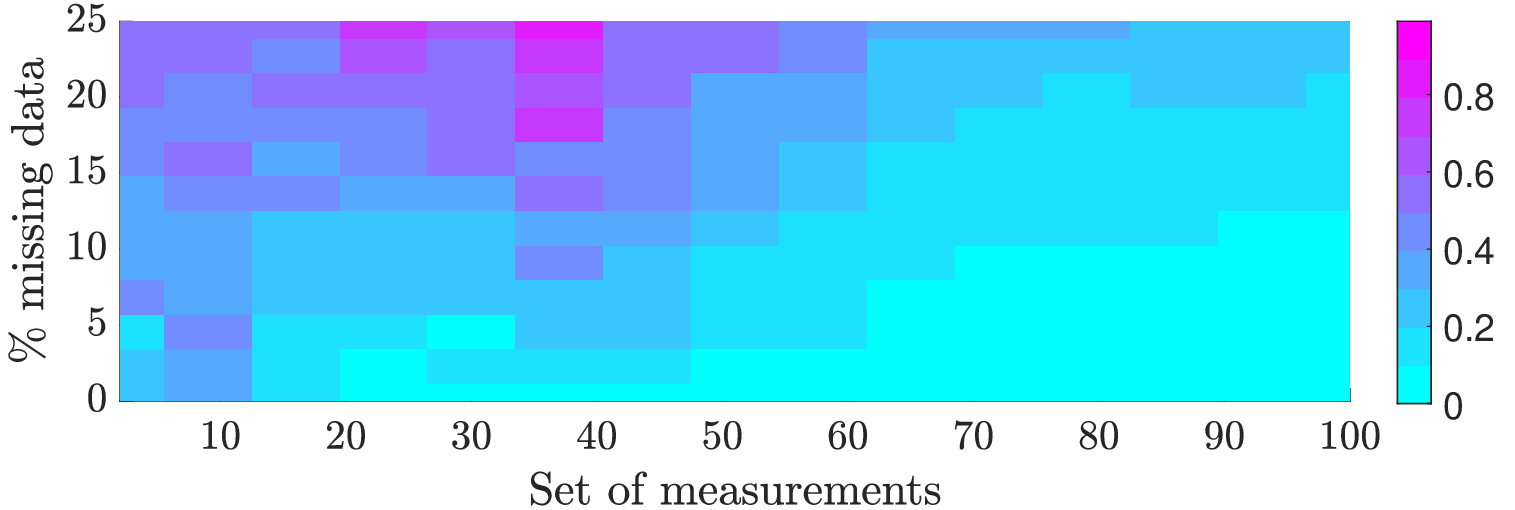}}
  \subfigure[]{\includegraphics[width=0.4\textwidth]{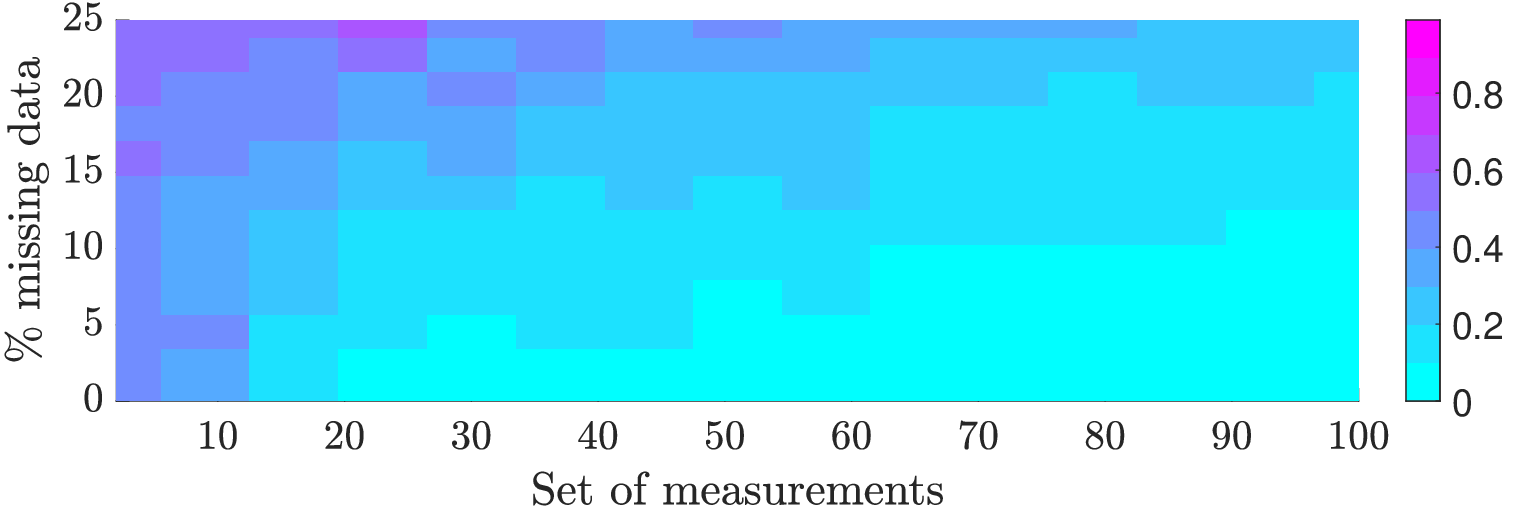}}
  \vspace{-.3cm}
  \caption{Case - 9-bus: Relative error (RE) over 50 trials for the estimation of the sensitivity matrix under a different number of measurements and different percentages of missing data: (a) RE for the least-square estimator, and (b) RE for the low-rank approach.}
  \vspace{-.25cm}
  \label{fig:Figure_Tra_missing}
\end{figure}

Further, in order to assess the performance of~\eqref{eq:noise_nuclear_outliers} in the case of  outliers, we replicated the previous case for the 9-bus transmission system but with random  outliers in the measurements. Figure \ref{fig:Figure_Tra_outliers} presents the results for the proposed method and the least-squares approach. Again, the proposed method outperforms the least-squares approach, and provides better estimates than the DC model-based method. In order to assess the performance of \eqref{eq:noise_nuclear_outliers_2} in the case of missing measurements, we replicated the case for the 9-bus transmission network. Figure \ref{fig:Figure_Tra_missing} presents the results for the proposed method and the least-squares approach where the LR method outperforms the LSE approach, in a case when different values of percentages of the data in $\Delta \mathbf{F}_k$ are missing.

%----------------------------------------
%
\subsection{Online Estimation}

\begin{figure}[!ht]
\centering
  \subfigure[]{\includegraphics[width=0.45\textwidth]{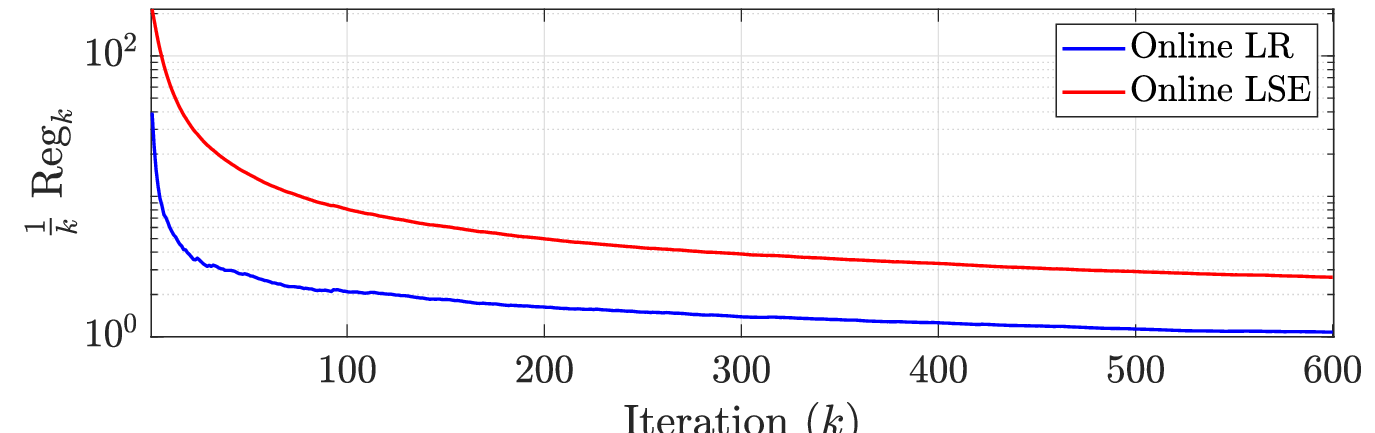}}
  \subfigure[]{\includegraphics[width=0.45\textwidth]{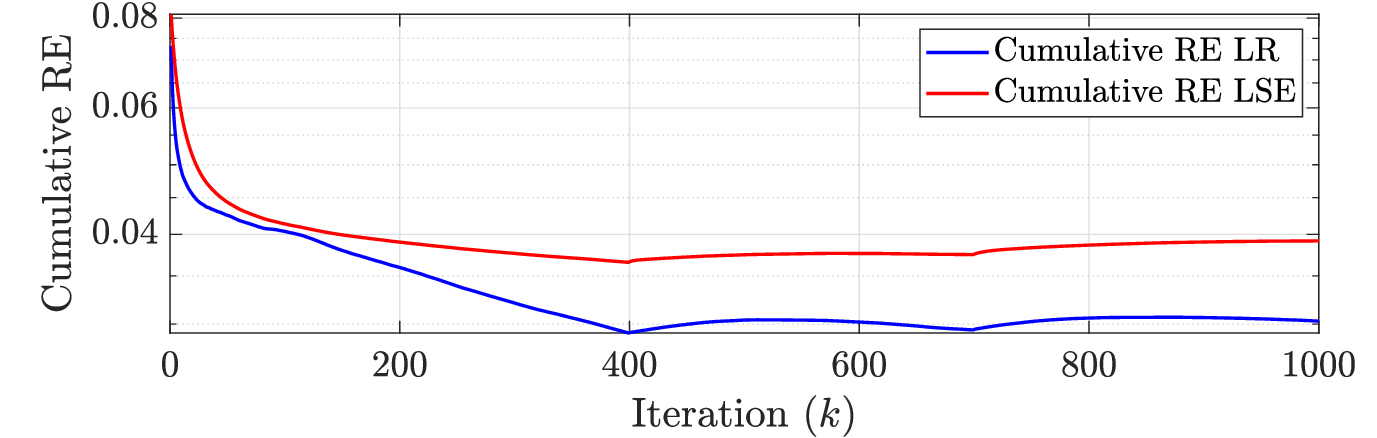}}
  \vspace{-.4cm}
  \caption{(a) Evolution $(1/k) \sum_{i=1}^k [f(\mathbf{x}_i) - f(\mathbf{x}^*_i)]$ for the online robust estimation of the sensitivity matrix in 9-bus transmission system using LR (low-rank) method and LSE (least-square estimation) approach. (b) Cumulative sum of the relative errors (RE) over $k$, i.e., $(1/k) \sum_{i=1}^k \mathrm{RE}_2$, for the online robust estimation of the sensitivity matrix in 9-bus transmission system, when there are changes of topology ($k = 400$ change reactance of line 5 and $k=700$ change reactance of line 8), using LR (low-rank) method and LSE (least-square estimation) approach.}
  \vspace{-.25cm}
  \label{fig:online_changetopo}
\end{figure}

As an example of an application of Algorithm~\ref{algo}, we consider an online robust estimation of the sensitivity matrix for the 9-bus transmission system. Relative to the test case presented in Section \ref{transmission_ex}, the nominal power injections at the nodes are now changing over time as in~\cite{dall2016optimal}.  Figure \ref{fig:online_changetopo}(a) shows the dynamic regret $(1/k) \text{Reg}_k$, when a window of 18 measurements is used. Based on Theorem~\ref{thm:convergence}, in the current setting the limiting behavior of $(1/k)\text{Reg}_k$ is $\mathcal{O}(1)$. Indeed, we can see that an asymptotic error is decreasing with the time index. Figure \ref{fig:online_changetopo}(b) presents the cumulative sum of the relative error (RE) over $k$, i.e., $(1/k) \sum_{i=1}^k \mathrm{RE}_2$, for the online robust estimation of the sensitivity matrix in the 9-bus transmission system, when there are changes of topology. In this case, we change the reactance of line 5 at $k = 400$, and the reactance of line 8 at $k=700$.

%------------------------------------------
% CONCLUSION
%------------------------------------------
\section{Conclusion}

This paper proposed a method to estimate sensitivities in a power grid by leveraging a nuclear norm minimization approach and sparsity-promoting regularization functions. The proposed methodology is applicable to the estimation of various sensitivities at the transmission level. Relative to a least-squares estimation method, the proposed approach allows to obtain meaningful estimates of the sensitivity matrix even when measurements are correlated. The method can identify outliers due to faulty sensors and is not deterred by missing measurements. An online proximal-gradient algorithm was proposed to estimate sensitivity matrices on-the-fly and enable operators to maintain up-to-date information on sensitivities under dynamic operating conditions. 

%------------------------------------------
% APPENDIX
%------------------------------------------
\appendix
\emph{Proof of Theorem~\ref{thm:convergence}}. 
Since $s_k$ has a $L_k$-Lipschitz continuous gradient, i.e., $L_k \geq \norm{\mathbf{A_{Ps}}_{,k}^\top \mathbf{A_{Ps}}_{,k}}$, then:
\vspace{-.2cm}
\begin{multline}
    s_k(\mathbf{x}_k) \leq s_k(\mathbf{x}_{k-1}) + \inp{\nabla_\mathbf{x}s_k(\mathbf{x}_{k-1})}{\mathbf{x}_k - \mathbf{x}_{k-1}}\\
    + \frac{L_k}{2}\norm{\mathbf{x}_k - \mathbf{x}_{k-1}}^2.
    \label{proof_1}
\end{multline}
Using the convexity of $s_k$ we also have that,
\vspace{-.2cm}
\begin{equation}
    s_k(\mathbf{x}_{k-1}) \leq s_k(\mathbf{x}_k^*) + \inp{\nabla_\mathbf{x}s_k(\mathbf{x}_{k-1})}{\mathbf{x}_{k-1} - \mathbf{x}_{k}^*}.
    \label{proof_2}
\end{equation}
Therefore, putting \eqref{proof_1} and \eqref{proof_2} together, one arrives at:
\vspace{-.2cm}
\begin{multline}
    s_k(\mathbf{x}_k) \leq s_k(\mathbf{x}_k^*) + \inp{\nabla_\mathbf{x}s_k(\mathbf{x}_{k-1})}{\mathbf{x}_{k} - \mathbf{x}_{k}^*}\\
    + \frac{L_k}{2}\norm{\mathbf{x}_k - \mathbf{x}_{k-1}}^2.
    \label{proof_3}
\end{multline}
% Definition of subgradient
On the other hand, for the non-differentiable function $g_k$, we can leverage \cite[Theorem 3.36]{Beck_prox}. Let $\phi_1(\mathbf{z})$, $\phi_2(\mathbf{z})$: $\mathbb{E} \rightarrow (-\infty, \infty]$ be proper convex functions, and let $\mathbf{z} \in \text{int}(\text{dom}(\phi_1)) \cap \text{int}(\text{dom}(\phi_2))$. For $\phi(\mathbf{z}) 
:= \phi_1(\mathbf{z})+\phi_2(\mathbf{z})$, then $\partial \phi(\mathbf{z}) = \partial \phi_1(\mathbf{z})+\partial \phi_2(\mathbf{z})$.
Based on \cite[Theorem 3.36]{Beck_prox}, we have that,
\vspace{-.2cm}
\begin{equation*}
     0 \in \partial \Big \{ g_k(\mathbf{x}_k) + \frac{1}{2 \alpha} \norm{\mathbf{x}_k - \mathbf{y}_k}^2 \Big \}, 
\end{equation*}
which implies:
\vspace{-.2cm}
\begin{equation*}
    -\bvarphi \in \partial g_k(\mathbf{x}_k), \qquad  \bvarphi \in \partial \Big \{\frac{1}{2 \alpha} \norm{\mathbf{x}_k - \mathbf{y}_k}^2\Big \}.
\end{equation*}
\noindent
Since $\bvarphi \in \partial \Big \{\frac{1}{2 \alpha} \norm{\mathbf{x}_k - \mathbf{y}_k}^2\Big \}$, the following holds:
\vspace{-.2cm}
% Replace gt by the definition of a subgradient
\begin{equation}
    \frac{1}{2 \alpha} \norm{\mathbf{x}_k - \mathbf{y}_k}^2 + \inp{\bvarphi}{\mathbf{q} - \mathbf{x}_k} - \frac{1}{2 \alpha} \norm{\mathbf{q} - \mathbf{y}_k}^2 \leq 0 \quad \forall \; \mathbf{q}.
    \label{sub_ineq}
\end{equation}
Furthermore, since \eqref{sub_ineq} holds for all $\mathbf{q}$, we can define $\mathbf{q} = \mathbf{y}_k + \alpha \bvarphi$. Then, \eqref{sub_ineq} can be written as
\vspace{-.15cm}
\begin{equation}
    \frac{1}{2 \alpha} (\norm{\mathbf{x}_k - \mathbf{y}_k}^2 + \alpha^2 \norm{\bvarphi}^2) - \inp{\bvarphi}{\mathbf{x}_k-\mathbf{y}_k} \leq 0.
    \label{include_q}
\end{equation}
Now using \eqref{sub_ineq} and \eqref{include_q} we get that, 
\vspace{-.15cm}
\begin{equation}
    \frac{\mathbf{y}_k - \mathbf{q}}{\alpha} \in \partial g_k(\mathbf{x}_k).
    \label{subgra_g}
\end{equation}
By using the subgradient defined in \eqref{subgra_g} and the update in \eqref{prox_update}, the following inequality for $g_k$ can be obtained,
\vspace{-.2cm}
\begin{multline}
    g_k(\mathbf{x}_k) \leq g_k(\mathbf{x}_k^*) - \frac{1}{\alpha} \inp{\mathbf{x}_{k-1}-\mathbf{q}}{\mathbf{x}_k^*-\mathbf{x}_k}\\
    + \inp{\nabla_\mathbf{x}s_k(\mathbf{x}_{k-1})}{\mathbf{x}_k^*-\mathbf{x}_k}.
    \label{proof_4}
\end{multline}

Adding \eqref{proof_3} and \eqref{proof_4} we obtained
\vspace{-.2cm}
\begin{multline*}
    s_k(\mathbf{x}_k) +  g_k(\mathbf{x}_k) \leq s_k(\mathbf{x}_k^*) + g_k(\mathbf{x}_k^*) + \frac{L_k}{2}\norm{\mathbf{x}_k - \mathbf{x}_{k-1}}^2\\ +\inp{\nabla_\mathbf{x}s_k(\mathbf{x}_{k-1})}{\mathbf{x}_{k}^* - \mathbf{x}_{k}} + \inp{\nabla_\mathbf{x}s_k(\mathbf{x}_{k-1})}{\mathbf{x}_k-\mathbf{x}_k^*}\\
    + \frac{1}{\alpha} \inp{\mathbf{x}_{k-1}-\mathbf{q}}{\mathbf{x}_k-\mathbf{x}_k^*},
\end{multline*}
and therefore,
\vspace{-.2cm}
\begin{multline}
    \hspace{-3mm} f_k(\mathbf{x}_k) \leq f_k(\mathbf{x}_k^*) + \left(  \frac{L_k}{2} - \frac{1}{\alpha} \right) \norm{\mathbf{x}_{k}-\mathbf{x}_k^*}^2 + 
    \\ \frac{L_k}{2} \norm{\mathbf{x}_{k-1}-\mathbf{x}_k^*}^2 + \left( \frac{1}{\alpha} - \frac{L_k}{2} \right) \inp{\mathbf{x}_{k-1} - \mathbf{q}}{\mathbf{x}_k - \mathbf{x}_{k}^*}.
    \label{proof_10}
\end{multline}
\noindent
Set $\alpha \leq \frac{1}{\text{max}\{L_k\}}$; in particular, let $\alpha = \frac{1}{\text{max}\{L_k\}}-w^2$, for $w \in \mathbb{R}$. Then:
\vspace{-.2cm}
\begin{equation}
    \frac{L_k}{2} - \frac{1}{\alpha} \leq \frac{\alpha \; \text{max}\{L_k\}-2}{2 \alpha} \leq -\frac{1}{2 \alpha}- \frac{w^2 \; \text{max}\{L_k\}}{2 \alpha}.
    \label{constant}
\end{equation}
% Using Chauchy-Schawarz inequality
Also notice that,
\vspace{-.2cm}
\begin{multline*}
    \left( \frac{1}{\alpha} - \frac{L_k}{2} \right) \inp{\mathbf{x}_{k-1} - \mathbf{q}}{\mathbf{x}_k - \mathbf{x}_{k}^*} \leq 
    \\ \left( \frac{1}{\alpha} - \frac{L_k}{2} \right) \norm{\mathbf{x}_{k-1} - \mathbf{q}}(\norm{\mathbf{x}_k - \mathbf{x}_{k}^*}+\omega_k).
\end{multline*}
Based on \eqref{sub_ineq} holds for all $\mathbf{q}$, let $R$ be the diameter of $\mathcal{X}$. Therefore,
\vspace{-.15cm}
\begin{equation}
    \left( \frac{1}{\alpha} - \frac{L_k}{2} \right) \norm{\mathbf{x}_{k-1} - \mathbf{q}}(\norm{\mathbf{x}_k - \mathbf{x}_{k}^*}+\omega_k) \leq \Phi R (R + \omega_k),
    \label{bound_3term}
\end{equation}
where $\Phi := \frac{1}{\alpha} - \frac{\text{min}\{L_k\}}{2}$. Then, based on \eqref{constant} and \eqref{bound_3term}, and neglecting constant terms, we can write \eqref{proof_10} as, 
\vspace{-.2cm}
\begin{multline}
    f_k(\mathbf{x}_k) - f_k(\mathbf{x}_k^*) \leq - \frac{1}{2 \alpha}\norm{\mathbf{x}_k - \mathbf{x}_k^*}^2 + \frac{1}{2 \alpha}\norm{\mathbf{x}_{k-1} - \mathbf{x}_k^*}^2 \\ + \Phi R (R + \omega_k).
    \label{proof_5}
\end{multline}

By adding and subtracting $\mathbf{x}^*_{k-1}$ in the last term on the right hand side of \eqref{proof_5}, and adding it from $i=1, \dots, k$, we can write the right hand side of the equation as:
\vspace{-.2cm}
\begin{equation}
    \sum_{i=1}^k \left\{ -\frac{1}{2 \alpha}\norm{\mathbf{x}_i - \mathbf{x}_i^*}^2 + \frac{1}{2 \alpha}\norm{\mathbf{x}_{i-1} - \mathbf{x}_{i-1}^* + \mathbf{x}_{i-1}^* - \mathbf{x}_i^*}^2 \right\},
    \label{proof_6}
\end{equation}
where the second term in \eqref{proof_6} can be expanded as:
\vspace{-.2cm}
\begin{multline}
    \sum_{i=1}^k \bigg\{ -\frac{1}{2 \alpha}\norm{\mathbf{x}_i - \mathbf{x}_i^*}^2 + \frac{1}{2 \alpha} \norm{\mathbf{x}_{i-1} - \mathbf{x}_{i-1}^*}^2 + \\ 
    \frac{1}{\alpha} \norm{\mathbf{x}_{i-1} - \mathbf{x}_{i-1}^*} \norm{\mathbf{x}_i^* - \mathbf{x}_{i-1}^*} + \frac{1}{2 \alpha} \norm{\mathbf{x}_i^* - \mathbf{x}_{i-1}^*}^2 \bigg\}.
    \label{proof_7}
\end{multline}
The first two terms in \eqref{proof_7} correspond a telescoping series; and, by using the definition of $\omega_k$, \eqref{proof_5} can be rewritten as follows:
\vspace{-.2cm}
\begin{multline}
    \sum_{i=1}^k [f_i(\mathbf{x}_i) - f_i(\mathbf{x}_i^*)] \leq  \frac{1}{2 \alpha} \norm{\mathbf{x}_{0} - \mathbf{x}_{0}^*}^2 -\frac{1}{2 \alpha}\norm{\mathbf{x}_k - \mathbf{x}_k^*}^2 \\ + \frac{1}{\alpha} \sum_{i=1}^k \omega_i (\norm{\mathbf{x}_{i-1} - \mathbf{x}_{i-1}^*} + \Phi R) + \frac{1}{2 \alpha} \sum_{i=1}^k \omega_i^2 + k \Phi R^2. \nonumber 
    \label{proof_8}
\end{multline}
Since $\mathcal{X}$ is compact, we can upper bound $\norm{\mathbf{x}_k - \mathbf{x}_k^*}$ by $R$, and thus the result follows.

\bibliographystyle{IEEEtran}
\bibliography{References}

\end{document}